\documentclass[letterpaper, 12 pt, draftclsnofoot, onecolumn]{ieeeconf}  

\IEEEoverridecommandlockouts                              
\usepackage[T1]{fontenc}
\usepackage[latin9]{inputenc}
\usepackage{amsmath}
\usepackage{amssymb}
\usepackage{graphicx}
\usepackage{dsfont}

\makeatletter

\pdfpageheight\paperheight
\pdfpagewidth\paperwidth

\@ifundefined{showcaptionsetup}{}{%
 \PassOptionsToPackage{caption=false}{subfig}}
\usepackage{subfig}
\makeatother

\newcommand{\Ni}{\mathcal{N}_i}
\newcommand{\co}[1]{\text{conv}\left\{#1\right\}}

\newcommand{\E}{\mathcal{E}}
\newcommand{\V}{\mathcal{V}}
\newcommand{\Tr}{\text{trace}}
\newcommand{\A}{\mathcal{A}}
\newcommand{\Adag}{\A^\dag}
\newcommand{\G}{\mathcal{G}}

\newcommand{\R}{\mathbb{R}}
\newcommand{\Sym}{\mathbb{S}}
\newcommand{\maximize}{\text{maximize}}
\newcommand{\IP}[1]{\left< #1 \right>}
\def\QED{~\rule[-1pt]{5pt}{5pt}\par\medskip}

\newtheorem{rem}{Remark}
\newtheorem{lem}{Lemma}

\newtheorem{thm}{Theorem}

\newtheorem{coro}{Corollary}
\newenvironment{IEEEproof}{{\bf Proof: }}{ \hfill \QED}
\newtheorem{example}{Example}

\newtheorem{assumption}{Assumption}

\begin{document}

\title{A convex approach to consensus on $SO(n)$}

\author{Nikolai Matni and Matanya B Horowitz
 \thanks{N. Matni is with the Department of Control and Dynamical Systems, California Institute of Technology, Pasadena, CA.
 \tt{\small nmatni@caltech.edu}.}
\thanks{M. Horowitz is with the Department of Control and Dynamical Systems, California Institute of Technology, Pasadena, CA.
 \tt{\small mhorowit@caltech.edu}.}
}
\maketitle
\begin{abstract}
This paper introduces several new algorithms for consensus over the special orthogonal group. By relying on a convex relaxation of the space of rotation matrices, consensus over rotation elements is reduced to solving a convex problem with a unique global solution. The consensus protocol is then implemented as a distributed optimization using (i) dual decomposition, and (ii) both semi and fully distributed variants of the alternating direction method of multipliers technique -- all with strong convergence guarantees. The convex relaxation is shown to be exact at all iterations of the dual decomposition based method, and exact once consensus is reached in the case of the alternating direction method of multipliers.  Further, analytic and/or efficient solutions are provided for each iteration of these distributed computation schemes, allowing consensus to be reached without any online optimization. Examples in satellite attitude alignment with up to 100 agents, an estimation problem from computer vision, and a rotation averaging problem on $SO(6)$ validate the approach.
\end{abstract}

\section{Introduction}
\label{sec:intro}
Optimization, coordination, and consensus over the group of rotation matrices (i.e. over elements of $SO(n)$) is a problem of fundamental importance in a wide range of applications, from satellite attitude and spin estimation \cite{SPW14,P09},  vehicle coordination \cite{BCM09}, frequency synchronization \cite{S00},  (distributed) visual pose estimation \cite{TVT08,TV09,horowitz2014convex} and protein folding \cite{LSS10}.

Traditional consensus in Euclidean space has a rich history \cite{OSFM07}, with results offering strong convergence guarantees under mild and realistic assumptions using \emph{purely local} protocols.  Further, with the resurfacing of distributed optimization methods such as dual decomposition \cite{PC06} and the alternating direction method of multipliers \cite{BoydADMM}, these methods have been successfully applied to achieve so-called ``fast consensus'' in a distributed \cite{EZDV11} and semi-distributed (i.e. sensor fusion) \cite{BoydADMM} setting.  Unfortunately, generalizing these approaches to manifolds, even those with as much structure as the group of rotation matrices, has proven non-trivial.

Fortunately, although the space of rotation matrices is highly non-linear, it is a Lie Group. This structure has allowed for sophisticated consensus methods to be created with so-called ``almost-global'' convergence -- that is to say the only \emph{stable} stationary points of the protocol are global minimizers of the underlying optimization.  These consensus schemes fall under two broad categories, namely \emph{intrinsic} and \emph{extrinsic} approaches.  

In the former, the consensus and/or optimization is performed directly on the manifold with respect to its geodesic distance: although many such methods exist in the literature\footnote{The interested reader should refer to the references of \cite{TAV12} for a more detailed overview of prior work.}, those found in \cite{TAV12} appear to be the most general and offer the strongest convergence guarantees.  This technique is able to achieve almost-global convergence, despite the inherent non-convexity of the underlying optimization, through an appropriate re-weighting of the consensus problem's objective function.  

The category of extrinsic protocols\footnote{Once again, a more exhaustive overview of the state of the art of extrinsic methods can be found in \cite{SS09}, and the references therein.} rely on embedding the underlying manifold into a Euclidean space, applying standard consensus protocols in this space, and then projecting the resulting states back onto the manifold \cite{SS09}.  These methods can also be constructed to provide almost-global convergence, but are highly dependent on the choice of embedding and projection methods, and do not easily generalize.

A common theme in all of the above methods is the use of sophisticated tools from differential geometry and Lie Group theory -- although elegant and appropriate for an applied mathematics community, such mathematical overhead may prove to be a barrier for the adoption of such methods by engineering practitioners.  

\textbf{Contributions: }This paper addresses all of these issues by exploiting a recent characterization of the convex hull of $SO(n)$ \cite{SPW14}.  In particular, we reduce the consensus optimization to a \emph{convex} one with a \emph{unique global minimizer} (under a mild technical assumption), and leverage distributed optimization techniques such as dual decomposition \cite{PC06} and the alternating direction method of multipliers \cite{BoydADMM} to implement a consensus protocol.  We further show that this global optimizer lies in $SO(n)$, rather than $\co{SO(n)}$, and thus the convex relaxation is indeed \emph{exact} once consensus is achieved.

We then prove that for a dual decomposition based method, the consensus path (i.e. all of the iterates of the protocol) are in fact in $SO(n)$, whereas for the ADMM protocol, we provide empirical evidence that an appropriate choice of algorithm parameters lead to iterates that are in or near $SO(n)$.  Finally, we show that our algorithms can be efficiently implemented, with individual agents only needing to perform an eigen-decomposition of a symmetric matrix, and/or a projection of a vector onto the probability simplex.


 We illustrate the usefulness of our approach on a rotation averaging problem (motivated by satellite formation control) for $N=8$, $N=50$ and $N=100$ agents arranged in a \emph{ring topology}, as well as on pose estimation problems based on observations given by several cameras, and on a consensus problem in $SO(6)$.   Finally, it is worth noting that if the results on the spectrahedral representation of the convex hull of $SO(n)$ from \cite{SPW14} and the well established distributed computation techniques we rely upon are taken as fact, the derivations require only elementary linear algebra and basic convex analysis.

\textbf{Paper outline:}  In Section \ref{sec:problem}, we begin by fixing notation and terminology, and then formulate the consensus problem over $SO(n)$ as a convex optimization over the convex hull of the group of rotations.  In particular, we show that the problem reduces to optimizing a linear functional defined by the average of the agents' data over this space.  In Section \ref{sec:prelims}, we provide an overview of the required results on the spectrahedral representation of the convex hull of $SO(n)$ \cite{SPW14}, and on distributed optimization techniques \cite{PC06,BoydADMM}.  In Section \ref{sec:consensus} we derive our consensus algorithms, and in Section \ref{sec:exp_results} we illustrate their effectiveness on several challenging consensus problems.  We end in Section \ref{sec:conclusion} with conclusions and directions for future work.

\section{Problem Formulation}
\label{sec:problem}
\subsection{Notation and Terminology}
We denote the Euclidean norm of a vector $x\in \R^n$ by $\|\cdot\|$.  The space of $n\times n$ symmetric matrices is denoted by $\Sym^n$, and is equipped with the natural inner product $\left<A,B\right>=\Tr (A^\top B) = \Tr (A B)$ which induces the Frobenius norm of an element $M\in \Sym^n$ as $\|M\|_F^2 = \Tr \left(M^2\right)$.  The set of $n$ dimensional rotation matrices, $SO(n)$, is defined as
\begin{equation}
SO(n):= \{ R \in \R^{n \times n} \ | \ R^\top R = I, \ \det R = 1\}.
\end{equation}

The convex hull operator is denoted by $\co{\cdot}$, and we use $M\succeq 0$ to indicate that $M\in\Sym^n$ is positive semi-definite.  As is standard, we refer to the set of unit trace positive semi definite matrices, that is to say those $M \in \Sym^n$ satisfying $M\succeq 0$, $\Tr M = 1$, as the \emph{free spectrahedron} -- this is the generalization of the probability simplex to $\Sym^n$.

\subsection{The Generalized Average}
Consider $N$ agents, each with state $R^i_t \in SO(n)$, interacting according to a graph $\G = (\V,\E)$, where $(i,j) \in \E \subset 2^\V$ if and only if agent $j$ can communicate its state to agent $i$ (note we assume $(i,i)\in\E$ for all $i\in\V$).  Further define the set of incoming neighbors at node $i$ as $\Ni := \{j \ | \ (i,j) \in \E\}$ to be those nodes that can communicate their state directly to node $i$ according to the graph topology.

In this paper, we consider the problem of implementing a (semi) distributed algorithm such that each state $R^i_t$ converges to a consensus value $\bar{R}$, which we define as the generalized average rotation of the system given by
\begin{equation}
\bar{R} := \arg \min_{R\in SO(n)} \sum_{i=1}^N \|C^i - RB^i\|_F^2,
\label{eq:non_convex}
\end{equation}
where each $C^i$ and $B^i$ are appropriately sized local data matrices.

Before justifying the use of this cost function and the generalized average terminology, we begin by providing a \emph{convex} reformulation of optimization \eqref{eq:non_convex}.

\begin{lem}
For $\bar{R}$ defined as in \eqref{eq:non_convex}, we have that
\begin{equation}
	\bar{R} = \arg \max_{R\in \co{SO(n)}} \IP{\frac{1}{N}\sum_{i=1}^N C^i(B^i)^\top,R}.
	\label{eq:convex}
\end{equation}
\end{lem}
\begin{IEEEproof}
Expanding the objective of optimization \eqref{eq:non_convex}, we may rewrite it as
\begin{equation}
\|C^i\|_F^2 - 2\left<C^i,RB^i\right> + \|RB^i\|_F^2.
\end{equation}
Noting that elements of $SO(n)$ leave the Frobenius norm invariant, the only non-constant term is the middle inner-product.  Multiplying by appropriate constants and taking the negative sign into account, we then obtain the following non-convex optimization
\begin{equation}
	\bar{R} = \arg \max_{R\in {SO(n)}} \IP{\frac{1}{N}\sum_{i=1}^N C^i(B^i)^\top,R}.
\end{equation}
However, noting that we are optimizing a linear functional over a basic semi-algebraic set, we may replace the feasible set with its convex hull, yielding the desired result.
\end{IEEEproof}

Thus from the objective of \eqref{eq:convex}, we see that the optimal solution is the element of $SO(n)$ that maximizes a linear functional defined by the average of a suitable function of the data -- hence the generalized average moniker.  We now present two examples to make these ideas concrete.

\begin{example}[Rotation Averaging]
Consider the case where $B^i = I$ and $C^i = R^i_0$, that is to say we wish to find the average value (with respect to the Frobenius norm) of the initial states of each node.  Then optimization \eqref{eq:convex} reduces to
\begin{equation}
\begin{array}{rcl}
\bar{R} 
&=& \arg \max_{R\in \co{SO(n)}}  \left<\frac{1}{N}\sum_{i=1}^N R^i_0,R\right> 
\end{array}
\end{equation}

\end{example}

 $\,$ 

\begin{example}[Multi-view pose estimation/localization]
Let $n=3$, and suppose that $B^i=X^i \in \R^{3\times m}$ are $m$ appropriately re-centered points obtained via either LIDAR or stereo camera observations, and that $C^i=Y^i \in \R^{3\times m}$ are the corresponding points of the internal model of the object being observed.  If all of the $Y^i$ are identical, this corresponds to using multiple cameras to estimate the pose of a single object, whereas different $Y^i$ could arise in, for example, a simultaneous localization and mapping (SLAM) scheme.  The optimization then becomes 
\begin{equation}
\begin{array}{rcl}
\bar{R} 
&=& \arg \max_{R\in \co{SO(n)}}  \left<\frac{1}{N}\left(\sum_{i=1}^N Y^i(X^i)^\top\right),R\right> 
\end{array}
\end{equation}
\end{example}

Based on the previous discussion, we let $D^i := C^i (B^i)^\top$ be a general data matrix available at node $i$, and consider solving the following optimization problem in a \emph{distributed manner} such that each local iterate $R^i_t$ converges to the desired generalized average
\begin{equation}
 \bar{R}:=\arg\max_{R\in \co{SO(n)}}  \sum_{i=1}^n \left<D^i,R\right>.
\label{eq:lin_opt}
\end{equation}

Our approach will be to exploit recent results on the spectrahedral representation of the \emph{convex hull} of $SO(n)$ \cite{SPW14}, and to leverage first order distributed optimization methods such as dual decomposition and alternating direction method of multipliers (ADMM) \cite{BoydADMM} to implement a consensus protocol that is guaranteed to converge to $\bar{R}$.  

In particular, we will consider two types of consensus settings: a \emph{completely distributed} setting in which information can only be exchanged between neighbors (i.e. node $i$ has access to $R^j$ if and only if $(i,j)\in\E$), and a \emph{semi-distributed} setting, in which information can not only be exchanged between neighbors, but can also be routed along so that each agent has access to every other agent's information.

\section{Preliminaries}
\label{sec:prelims}
\subsection{The convex hull of $SO(n)$}
In the previous section, we showed that our consensus problem amounts to optimizing a linear functional over the convex hull of $SO(n)$. However, unless this latter object has a tractable representation, we have not made much progress.  Fortunately, it has recently been shown that $\co{SO(n)}$ admits a \emph{spectrahedral} representation of dimension $2^{n-1} \times 2^{n-1}$. Although this scaling is poor in $n$, most applications have low $n$, e.g.,  rigid body transformations have $n=2,3$.

To that end, we recall the characterization \eqref{eq:conSOn} of $\co{SO(n)}$ from Corollary 1.6 of $\cite{SPW14}$.  In particular, for $1\leq i,j \leq n$, we have that
\begin{equation}
A_{ij} = -P^\top_{\text{even}}\lambda_i\rho_jP_{\text{even}},
\end{equation}
where $\lambda_i$ and $\rho_j$ are given by
\begin{equation}
\begin{array}{rcl}
\lambda_i &=& \overbrace{\begin{bmatrix} 1 & 0 \\ 0 & -1 \end{bmatrix} \otimes \cdots \otimes \begin{bmatrix} 1 & 0 \\ 0 & -1 \end{bmatrix} }^{i-1} \otimes \begin{bmatrix} 0 & -1 \\ 1 & 0 \end{bmatrix} \otimes \overbrace{\begin{bmatrix} 1 & 0 \\ 0 & 1 \end{bmatrix} \otimes \cdots \otimes \begin{bmatrix} 1 & 0 \\ 0 & 1 \end{bmatrix} }^{n-i} \\ \\
\rho_i &=& \underbrace{\begin{bmatrix} 1 & 0 \\ 0 & 1 \end{bmatrix} \otimes \cdots \otimes \begin{bmatrix} 1 & 0 \\ 0 & 1 \end{bmatrix} }_{i-1} \otimes \begin{bmatrix} 0 & -1 \\ 1 & 0 \end{bmatrix} \otimes \underbrace{\begin{bmatrix} 1 & 0 \\ 0 & -1 \end{bmatrix} \otimes \cdots \otimes \begin{bmatrix} 1 & 0 \\ 0 & -1 \end{bmatrix} }_{n-i} 
\end{array},
\label{eq:lamrho}
\end{equation}
and
\begin{equation}
P_{\text{even}} = \frac{1}{2} \begin{bmatrix} I + \left[\begin{smallmatrix} 1 & 0 \\ 0 & -1 \end{smallmatrix}\right] \otimes \cdots \otimes \left[\begin{smallmatrix} 1 & 0 \\ 0 & -1 \end{smallmatrix}\right] \\
\\
 I - \left[\begin{smallmatrix} 1 & 0 \\ 0 & -1 \end{smallmatrix}\right] \otimes \cdots \otimes \left[\begin{smallmatrix} 1 & 0 \\ 0 & -1 \end{smallmatrix}\right]
 \end{bmatrix}.
\end{equation}

The convex hull can then be written as
\begin{equation}
\co{SO(n)} = \left\{ \begin{bmatrix} \IP{A_{11},Z} & \IP{A_{12},Z} & \dots & \IP{A_{1n},Z} \\
\IP{A_{21},Z} & \IP{A_{22},Z} & \dots & \IP{A_{2n},Z} \\
\vdots & \vdots & \ddots &\vdots \\
\IP{A_{n1},Z} & \IP{A_{n2},Z} & \dots & \IP{A_{nn},Z} \end{bmatrix} \, | \, Z \succeq 0, \, \text{trace}Z=1 \right\}.
\label{eq:conSOn}
\end{equation}

For our purposes, however, it will suffice to note that $\co{SO(n)}$ can be written as
\begin{equation}
\co{SO(n)} = \left\{ \A(Z) \ | \ Z \succeq 0, \ \text{trace}(Z) =1 \right\},
\label{eq:simple_SOn}
\end{equation}
for the affine map $\A:\Sym^{2^{n-1}}\to \R^{n \times n}$ defined in terms of the matrices $A_{ij}$ as given in \eqref{eq:conSOn}.  Further let $\Adag:\R^{n \times n}\to\Sym^{2^{n-1}}$ denote its adjoint operator with respect to $\IP{\cdot,\cdot}$.  Patient computation reveals that 
\begin{equation}
\Adag(Y) = \sum_{i,j} A_{ij}Y_{ij}.
\label{eq:adjoint}
\end{equation}

\subsubsection{An eigenvalue problem}
We now exploit this characterization of $\co{SO(n)}$ to provide an analytic solution to the following optimization problem, which will be an essential element of our consensus algorithms.

\begin{lem}
\label{lem:analytic}
The optimal solution $R^*\in SO(n)$ to the optimization
\begin{equation}
\max_{R\in \co{SO(n)}} \left<D,R\right>
\label{eq:lem_opt}
\end{equation}
is given by
\begin{equation}
R^* = \A(\mu \mu^\top),
\end{equation}
where $\mu$ is the orthonormal eigenvector corresponding to the largest eigenvalue $\nu$ of $\Adag(D)$.
\end{lem}
\begin{IEEEproof}
Using the characterization of $\co{SO(n)}$ \eqref{eq:simple_SOn}, optimization \eqref{eq:lem_opt} is equivalent to
\begin{equation}
\max_{Z \succeq 0, \, \Tr Z = 1} \left<\Adag(D),Z\right>.
\end{equation}

It is then straightforward to verify that (i) Slater's condition is satisfied, and thus strong duality holds, and that (ii) $(\mu\mu^\top, \nu_{\max})$ form a Karush-Kuhn-Tucker pair for the optimization.  Finally, as $Z^*:=\mu\mu^\top$ is of rank 1, $\mathcal{A}(Z^*)$ is an extreme point of $\co{SO(n)}$, and thus an element of $SO(n)$.
\end{IEEEproof}

\begin{coro}
Suppose that the maximal eigenvalue $\nu$ of $\Adag(D)$ has multiplicity 1.  Then the unique optimal solution $Z^*$ to the optimization 
\begin{equation}
\max_{Z \succeq 0, \, \Tr Z = 1} \left<\Adag(D),Z\right>
\end{equation}
is given by $Z=\mu\mu^\top$.
\label{coro:uniqueness}
\end{coro}
\begin{IEEEproof}
As we are optimizing a linear functional over the free spectrahedron, the optimal solution must be an extreme point satisfying $Z^* = zz^\top$, $\|z\|=1$.  Thus, to show uniqueness it suffices to note that $\IP{zz^\top,\Adag(D)}=z^\top\Adag(D)z<\nu$ for all $z \neq \mu$ satisfying $\|z\| =1$.  
\end{IEEEproof}

\begin{rem}
The hypothesis that the maximal eigenvalue $\nu$ of $\Adag(D)$ has multiplicity 1 will be satisfied with probability one for randomly distributed data matrices.
\end{rem}


\subsection{Dual Decomposition}

Dual decomposition is a general technique to solve optimization problems
with separable objectives, and has the advantage of yielding a distributed algorithm. Specifically,
given the problem
\begin{eqnarray}
\min_{x_{i}\in\mathbb{R}^{n_{i}}} &  & f(x)=\sum_{i=1}^{N}f^{i}(x^{i})\label{eq:dual_decomp_original}\\
\text{s.t.} &  & Ax=b\nonumber 
\end{eqnarray}
where the optimization variable $x$ is partitioned according to the
objective as $x=\left[x_{1}^{T},\ldots,x_{N}^{T}\right]^{T}$, and the
constraint matrix $A$ can be also be partitioned according to this
structure, with $A=\left[A_{1},\ldots,A_{N}\right]$. The Lagrangian
for (\ref{eq:dual_decomp_original}) is then

\begin{eqnarray*}
L(x,\lambda) & = & \sum_{i=1}^{N}f^{i}(x^{i})+\lambda^{T}A_{i}x^{i}-\frac{1}{N}\lambda^{T}b\\
 & = & \sum_{i=1}^{N}L^{i}(x^{i},\lambda)
\end{eqnarray*}
with Lagrange multipliers introduced as $\lambda$. This gives rise
to the dual function
\[
g(\lambda)=\inf_{x}L(x,\lambda).
\]

This dual function may be optimized via gradient ascent. Due to the
separability of the Lagrangian, this yields the following procedure

\begin{eqnarray*}
x^{i}_{t+1} & = & \text{argmin}_{x^{i}}L^{i}(x^{i},\lambda_{t})\\
\lambda_{t+1} & = & \lambda_{t}+\alpha_{t}\left(Ax_{t+1}-b\right)
\end{eqnarray*}
where it is seen that the optimization may be performed in parallel
over the primal variables $x_{i}$ at each iteration $k$. At each
iteration, these updated primal variables are gathered for the update
of the dual variable, which may then be distributed to the computation
units. These dual variable updates may also be performed in parallel according to the separability of $A$. The result is a distributed algorithm for solving (\ref{eq:dual_decomp_original}), which can be shown to converge under certain conditions (which our problem will be shown to satisfy).


\subsection{Alternating Direction Method of Multipliers}
The Alternating Direction Method of Multipliers (ADMM) \cite{BoydADMM, eckstein1989splitting} provides a principled method for parallelization of convex problems. It is adopted here due to its particularly strong and general convergence guarantees, allowing us to inherit these guarantees in our consensus algorithm.

ADMM is a ``meta''-optimization scheme, where each step is carried out by solving a
convex optimization problem. Consider the optimization
\begin{equation}\label{eq:admm_standard_form}
	\begin{array}{ll}
		\mbox{minimize}   & f(x)+g(z) \\
		\mbox{subject to} & Ax+Bz=c
	\end{array}
\end{equation}
over the variables $x$ and $z$ and convex functions $f$ and $g$. Define
an augmented Lagrangian
\[
	L_{\alpha}=f(x)+g(z)+y^{T}\left(Ax+Bz-c\right)
		+\frac{\alpha}{2}\left\|Ax+Bz-c\right\|_{2}^{2},
\]
where $\alpha>0$ is an algorithm parameter, and $y$ is the dual variable
associated with the equality constraint. The constrained optimization is solved
through alternately minimizing the augmented Lagrangian over the primal
variables $x$, $z$, and updating the dual variable $y$,
\begin{eqnarray*}
x_{t+1} & := & \text{argmin}_{x}L_{\alpha}(x,z_{t},y_{t})\\
z_{t+1} & := & \text{argmin}_{z}L_{\alpha}(x_{t+1},z,y_{t})\\
y_{t+1} & := & y_{t}+\alpha\left(Ax_{t+1}+Bz_{t+1}-c\right).
\end{eqnarray*}

We will be concerned with the following two assumptions:
\begin{assumption}\label{as:admm1}
The (extended real valued) functions $f:\mathbb{R}^n \to \mathbb{R} \cup \{+\infty\}$ and $g : \mathbb{R}^m \to \mathbb{R} \cup \{+\infty\}$ are closed, proper, and convex.
\end{assumption}
\begin{assumption}\label{as:admm2}
The unaugmented Lagrangian has a saddle point.
\end{assumption}
If it can be demonstrated that the optimization problem obeys these assumptions, then the following general theorem is available:
\begin{thm}\label{thm:admm_convergence}
(See \cite{BoydADMM}) Given Assumptions \ref{as:admm1}, \ref{as:admm2} then the ADMM iterates satisfy the following:
\begin{itemize}
\item \textbf{Residual convergence}: $r_t \to 0$ as $t \to \infty$, i.e. the iterates approach feasibility
\item \textbf{Objective convergence}: $f(x_t) + g(z_t) \to p^*$ as $t\to\infty$, i.e. the objective function of the iterates approaches the optimal value
\item \textbf{Dual variable convergence}: $y_t \to y^*$ as $t \to \infty$, where $y^*$ is a dual optimal point
\end{itemize}
\end{thm}

\section{Consensus Algorithms}
\label{sec:consensus}
\subsection{General Approach}
As is standard in the use of distributed optimization algorithms for consensus, we will rewrite optimization \eqref{eq:lin_opt} in a manner more amenable to our purposes by introducing local variables at each node, and enforcing consistency according to the connectivity of the graph.  In particular, we now look to solve the following optimization in a distributed manner:
\begin{equation}
\begin{array}{rl}
\underset{R^i, \, i \in \V}{\maximize} & \sum_{i=1}^N \IP{D^i, R^i} \\
\text{s.t.} & R^i \in \co{SO(n)} \ \forall i \in \V \\
& R^i = R^j \ \forall (i,j) \in \E,
\end{array}
\label{eq:dist_op1t}
\end{equation}
or in its parameterized form,
\begin{equation}
\begin{array}{rl}
\underset{Z^i, \, i\in\V}{\maximize} & \sum_{i=1}^N \IP{\Adag(D^i), Z^i} \\
\text{s.t.} & Z^i \succeq 0, \ \Tr Z^i = 1 \ \forall i \in \V \\
& Z^i = Z^j \ \forall (i,j) \in \E.
\end{array}
\label{eq:dist_opt}
\end{equation}

It is evident that so long as the graph is strongly connected, that the solution to this optimization is identical to that of \eqref{eq:lin_opt}, and thus we work with this distributed optimization for the remainder of the paper.

\subsection{Completely Distributed Consensus via Dual Decomposition}
\label{sec:dual_decomposition}

We write the Lagrangian of optimization \eqref{eq:dist_opt} as
\begin{equation}
L(Z^i, Y^{ij}) = \sum_{i=1}^N \IP{\Adag(D^i), Z^i} - \sum_{(i,j) \in \E} \IP{Y^{ij},Z^i - Z^j}
\end{equation}
where each $Z^i$ is constrained to be positive semi-definite and of unit trace, and the $Y^{ij}$ are symmetric.  The gradient of the dual objective function with respect to $Y^{ij}$ is then easily verified to be $Z^j - Z^i$, allowing us to write the following dual-ascent algorithm for the distributed optimization\footnote{Note that our presentation is slightly non-standard as our primal problem is a maximization, and our dual a minimization.}:
\begin{equation}
\begin{array}{rcl}
Z^i_{t+1} &=& \underset{Z^i\succeq 0, \Tr Z^i =1}{\text{argmax}} \IP{\Adag\left(D^i\right) - \sum_{j:(i,j)\in\E} Y^{ij}_t,Z^i} \\
Y^{ij}_{t+1} &=& Y^{ij}_t - \alpha_t \left( Z^j_{t+1} - Z^i_{t+1}\right)
\end{array}
\label{eq:dualdecomp}
\end{equation}

Leveraging Lemma \ref{lem:analytic}, we can solve for the $Z^i$ update analytically through an eigvenvalue decomposition of a $2^{n-1}\times 2^{n-1}$ symmetric matrix, making this method computationally appealing for small to moderate $n$.

\subsubsection{Convergence Guarantees}
Assuming that the hypothesis of Corollary \ref{coro:uniqueness} holds at each iteration (which will be true for generic data), then the iterates are guaranteed to converge to the global optimizer $Z^*$ of \eqref{eq:dist_opt} (c.f. Prop. 3.9 of \cite{BT89}) for appropriately chosen step size $\alpha_t$.

\subsubsection{Iterate properties}
As the objective of the $Z^i_{t+1}$ optimization step in \eqref{eq:dualdecomp} is linear in $Z^i$, by Lemma \ref{lem:analytic}, each iterate is an element of $SO(n)$. Thus, this approach is well suited for use in ``true consensus'' protocols in which each iterate is required to be a valid rotation element.

\subsection{Completely Distributed Consensus via ADMM}
\label{sec:dist_admm}

Following the approach of \cite{EZDV11}, we introduce additional auxiliary variables $X^j$ in optimization \eqref{eq:dist_opt} and rewrite it as
\begin{equation}
\begin{array}{rl}
\underset{Z^i,X^j,\, i,j\in\V}{\maximize} & \sum_{i=1}^N \IP{\Adag(D^i), Z^i} \\
\text{s.t.} 	& Z^i \succeq 0, \ \Tr Z^i = 1 \ \forall i \in \V \\
		& X^j \succeq 0, \ \Tr X^j = 1 \ \forall j \in \V \\
		& Z^i = X^j \ \forall (i,j) \in \E.
\end{array}
\label{eq:distributed}
\end{equation}

We then write the augmented Lagrangian of optimization \eqref{eq:distributed} as
\begin{multline}
L(Z^i, X^j, Y^{ij}) = \sum_{i=1}^N \IP{\Adag(D^i), Z^i} - \dots \\ \hfill{} \sum_{(i,j) \in \E} \left[ \IP{Y^{ij},Z^i - X^j} + \frac{\alpha}{2}\|Z^i - X^j\|^2_F\right]
\end{multline}
where each $Z^i$ and $X^j$ are constrained to be positive semi-definite and of unit trace, and the $Y^{ij}$ are symmetric.  Exploiting the identity $\sum_{(i,j)\in\E} = \sum_{i=1}^N \sum_{j:(i,j)\in\E} = \sum_{j=1}\sum_{i:(j,i) \in \E}$, the ADMM algorithm then becomes
\begin{equation}
\begin{array}{rcl}
Z^i_{t+1} &=& \underset{Z^i\succeq 0, \Tr Z^i =1}{\text{argmax}} \IP{M^i_t,Z^i} - \frac{|\Ni|\alpha}{2}\|Z^i\|_F^2 \\
X^j_{t+1} &=& \underset{X^j\succeq 0, \Tr X^j =1}{\text{argmax}} \IP{N^j_{t},X^j} - \frac{|\mathcal{N}_j|\alpha}{2}\|X^j\|_F^2 \\
Y^{ij}_{t+1} &=& Y^{ij}_t - \alpha \left( Z^j_{t+1} - Z^i_{t+1}\right)
\end{array}
\label{eq:dist_admm}
\end{equation}
where
\begin{equation}
\begin{array}{rcl}
M^i_t &=& \Adag\left(D^i\right) - \sum_{j:(i,j)\in\E} \left[Y^{ij}_t-\alpha X^j_t\right] \\
N^j_t &=& \sum_{i:(j,i)\in\E} \left[Y^{ij}_t + \alpha Z^i_{t+1}\right].
\end{array}
\end{equation}

\subsection{Semi-Distributed Consensus via ADMM}
\label{sec:admm}

We first modify optimization \eqref{eq:dist_opt} by introducing a fusion variable $Z^0$:
\begin{equation}
\begin{array}{rl}
\underset{Z^0,Z^i,\, i\in\V}{\maximize} & \sum_{i=1}^N \IP{\Adag(D^i), Z^i} \\
\text{s.t.} & Z^i \succeq 0, \ \Tr Z^i = 1 \ \forall i \in \V\cup\{0\}\\
& Z^i = Z^0\ \forall i \in \V.
\end{array}
\label{eq:fusion_opt}
\end{equation}

We write the augmented Lagrangian of optimization \eqref{eq:dist_opt} as
\begin{multline}
L_\alpha(Z^i, Z, Y^{i}) = \sum_{i=1}^N \IP{\Adag(D^i), Z^i} - \dots \\ \sum_{i=1}^N \left(\IP{Y^{i},Z^i - Z^0} + \frac{\alpha}{2}\|Z^i - Z^0\|_F^2\right),
\end{multline}
where once again each $Z^i$ is constrained to be positive semi-definite and of unit trace, and the $Y^{i}$ are symmetric.  The ADMM algorithm then becomes
\begin{equation}
\begin{array}{rcl}
Z^i_{t+1} &=& \underset{Z^i\succeq 0, \Tr Z^i =1}{\text{argmax}} \IP{\Adag\left(D^i\right) - Y^{i}_t+\alpha Z^0_t,Z^i} \dots \\ && \hfill{} - \frac{\alpha}{2}\|Z^i\|_F^2 \\
Z^0_{t+1} &=& \underset{Z^0\succeq 0, \Tr Z^0 =1}{\text{argmax}}  \IP{\sum_{i=1}^N\left(Y^i_{t} + \alpha Z^i_{t+1}\right),Z^0} \dots \\ && \hfill{} - \frac{\alpha}{2}\|Z^0\|_F^2 \\
Y^{i}_{t+1} &=& Y^{i}_t - \alpha \left( Z^0_{t+1} - Z^i_{t+1}\right).
\end{array}
\label{eq:admm}
\end{equation}

\subsection{Convergence Guarantees and Efficient Updates}
Assumption 1 is satisfied as the objectives are linear, and Assumption 2 is satisfied as strong duality is easily verified to hold: thus this algorithm is guaranteed to converge for appropriately chosen regularization parameter $\alpha$ -- further as the global optimizer is unique and the objective functions continuous, convergence to an optimal cost implies convergence of the iterates to the optimal solution.

We note that as written, the $Z^i$ and $X^j$ updates do not admit an obvious efficient solution.  The next subsection will show how this problem can be reduced to a projection onto the probability simplex, which admits efficient algorithms that are of linear complexity in the dimension of the vector being projected \cite{MDP89}.

\subsubsection{$Z^i$ and $X^j$ updates}

We note that the $Z^i$ and $X^j$ updates in equations \eqref{eq:dist_admm} and \eqref{eq:admm} are all of the form
\begin{equation}
 \underset{W\succeq 0, \Tr W =1}{\text{argmax}}  \IP{\kappa\alpha T,W} - \frac{\kappa\alpha}{2}\|W\|_F^2 \\
 \label{eq:optimization}
\end{equation}
for appropriate constant matrix $T$ and scalar $\kappa$.  We begin by observing that the optimal solution $W^*$ to the above program is the projection of $T$ onto the free spectrahedron.

\begin{lem}
The optimal solution $W^*$ to optimization \eqref{eq:optimization} is given by the projection of $T$ onto the free spectrahedron.
\end{lem}
\begin{IEEEproof}
Multiplying the objective by $-\frac{2}{\kappa\alpha}$ and accordingly replacing the maximization with a minimization, we obtain
\begin{equation}
W^* = \underset{W\succeq 0, \Tr W =1}{\text{argmin}}-2\IP{T,W} + \|W\|_F^2.
\end{equation}
It suffices to note that adding $\|T\|^2_F$ to the objective does not change the optimal update's value, thus reducing the optimization to
\begin{equation}
W^* = \underset{W\succeq 0, \Tr W =1}{\text{argmin}}\|T - W\|_F^2,
\end{equation}
which is none other than the projection of $T$ onto the free spectrahedron.
\end{IEEEproof}

From this Lemma, the following two corollaries are immediate.

\begin{coro}
The $Z^i_{t+1}$ update as given in equation \eqref{eq:admm} is given by the projection of $Q^i_t$ onto the free spectrahedron, 
where
\begin{equation}
Q^i_t := \begin{cases} 
\sum_{i=1}^N \frac{1}{\alpha}Y^i_t + Z^i_{t+1} & \text{for $i=0$} \\
\frac{1}{\alpha}\left(\Adag(D^i) - Y^i_t\right) + Z^0_t & \text{for $i=1,\dots,N$}
\end{cases}
\end{equation}
\end{coro}

\begin{coro}
The $Z^i_{t+1}$ and $X^j_{t+1}$ update as given in equation \eqref{eq:dist_admm} are given by the projection of $\frac{1}{|\Ni|\alpha}M^i_t$ and $\frac{1}{|\mathcal{N}_j|\alpha}N^j_t$, respectively, onto the free spectrahedron.
\end{coro}

Finally, since in all cases, $T$ is symmetric, the aforementioned projection can be further reduced to a projection of the eigenvalues of $T$ onto the probability simplex.
\begin{lem}
Let $T = U^i_t \Lambda^i_t (U^i_t)^\top$ be a diagonalization of $T$.  Then the projection of $T$ onto the free spectrahedron is given by $U^i_t [\Lambda^i_t]^+ (U^i_t)^\top$, where $[\cdot]^+$ denotes a projection of the diagonal elements of a matrix onto the probability simplex.
\end{lem}
\begin{IEEEproof}
Follows immediately from diagonalizing $T$ and the invariance of the Frobenius norm under orthogonal transformations.
\end{IEEEproof}

\subsubsection{Iterate properties}
Although we are able to guarantee convergence to an element of $SO(n)$, the iterates may lie in the interior of $\co{SO(n)}$ -- the iterate optimizations for $Z^i_{t+1}$ and $X^j_{t+1}$ have a quadratic objective, and thus extreme point solutions are not guaranteed.  However, these ADMM schemes are still well suited for distributed computation on the group of rotations, in which only the consensus point is actually used.

In practice however, we have noticed that through an appropriate choice of $\alpha$, all iterates are in, or very near to, $SO(n)$.  In particular, by choosing a sufficiently low value of $\alpha$ (but still large enough to maintain good convergence), the objectives are dominated by their linear term, thus leading to (near) boundary solutions -- current work is aimed at formalizing this intuition.  Indeed this heuristic was used in the following examples to generate valid consensus paths within $SO(n)$.

\section{Experimental Results}
\label{sec:exp_results}
We tested our methods on a number of estimation and consensus problems. The dual decomposition and ADMM methods of \S\ref{sec:dual_decomposition} and \S\ref{sec:admm} were employed to synchronize the attitude for a formation of spacecraft. This is then followed by a pose estimation problem from computer vision, where the semi-distributed algorithm of \S\ref{sec:admm} allows for a consistent estimate of an object's orientation to be achieved from a collection of disparate measurements.  We finally end with an application of the dual decomposition method to a rotation averaging problem on $SO(6)$ to illustrate that our methods generalize seamlessly to higher dimensions.

\subsection{Satellite Attitude Synchronization via Dual Decomposition}
\begin{figure}
\begin{center}
\includegraphics[scale=1.4]{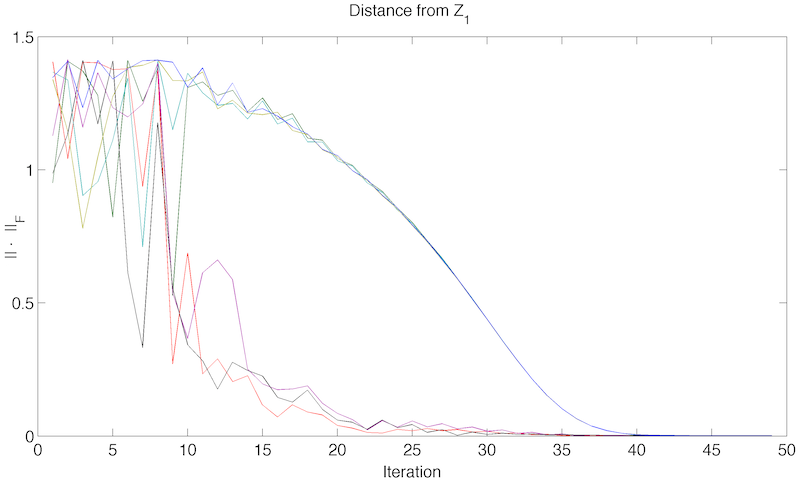}
\caption{\label{fig:dual_convergence}
Error between average rotation estimates for a satellite formation of $N=8$ satellites with a ring communication topology. The gradient descent weight is set as $\alpha=5$.
}
\end{center}
\end{figure}

Similar to \cite{TAV12} we propose to synchronize the attitude of a $N=8$ satellite formation. Each satellite is initialized with a random orientation, with the task being to find the average rotation, i.e. Example 1 of \S\ref{sec:problem}. The communication topology of the satellites is limited to a ring, and the problem solved via dual decomposition. The results are demonstrated in Figure \ref{fig:dual_convergence}, where convergent behavior is indeed observed. Computational time for each iteration for all satellites is typically under $1ms$ for unoptimized, interpreted (Matlab) code.

\begin{figure}
\begin{center}
\includegraphics[scale=.4]{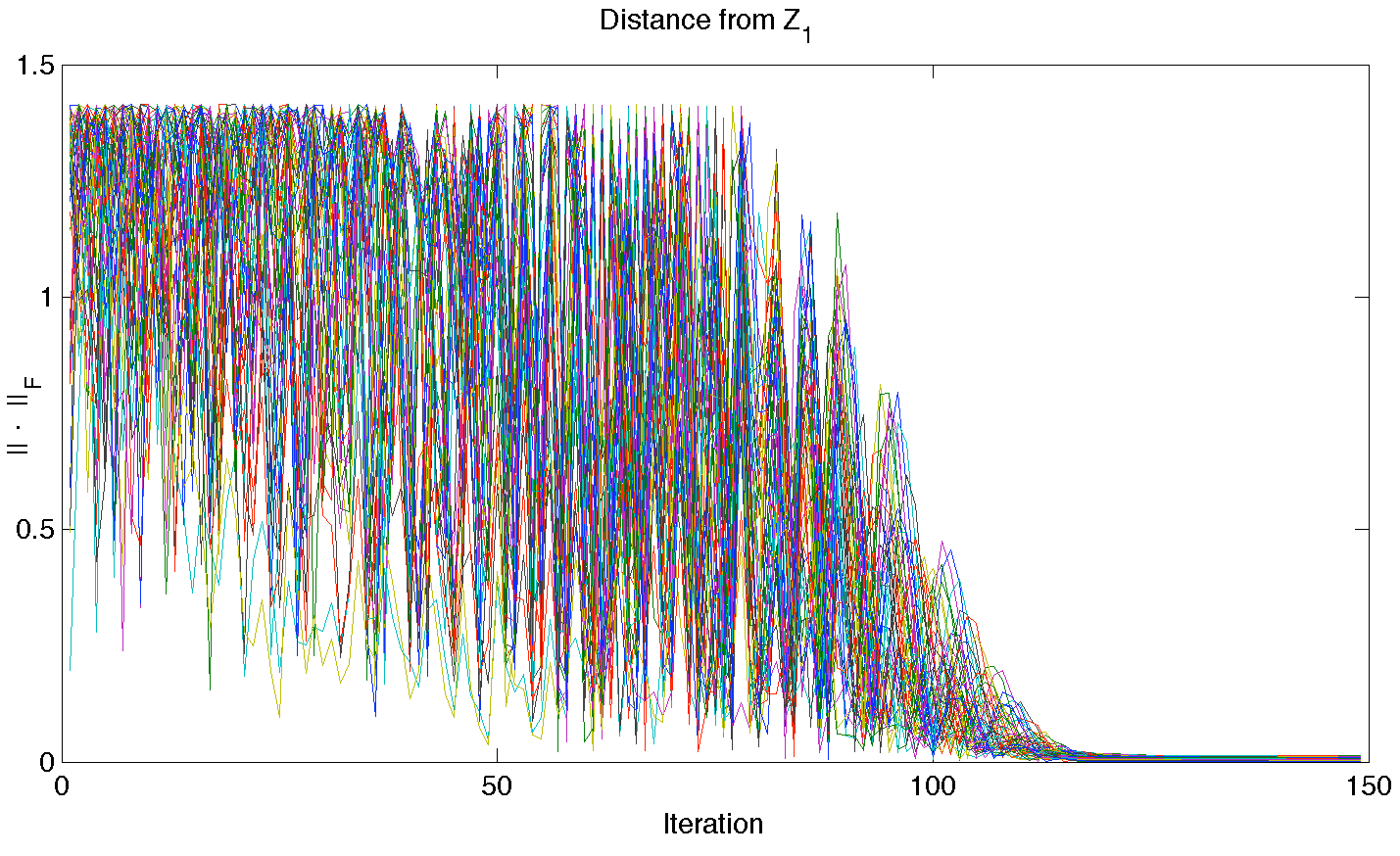}
\caption{\label{fig:dual_convergence_50}
Error between average rotation estimates for a satellite formation of $N=50$ with nearest neighbor and second nearest neighbor communication. The gradient weight is set as $\alpha=5$.
}
\end{center}
\end{figure}

A larger example is then tested, where $n=50$ and the communication topology is augmented to allow communication between both neighbors and second-neighbors. The resulting error trace is shown in Figure \ref{fig:dual_convergence_50}, where convergence is again observed. In practice, convergence is obtained for all examples for a wide range values $\alpha \in (.01 - 20)$.

\subsection{Attitude Synchronization via Distributed ADMM}
\begin{figure}
\begin{center}
\includegraphics[scale=.4]{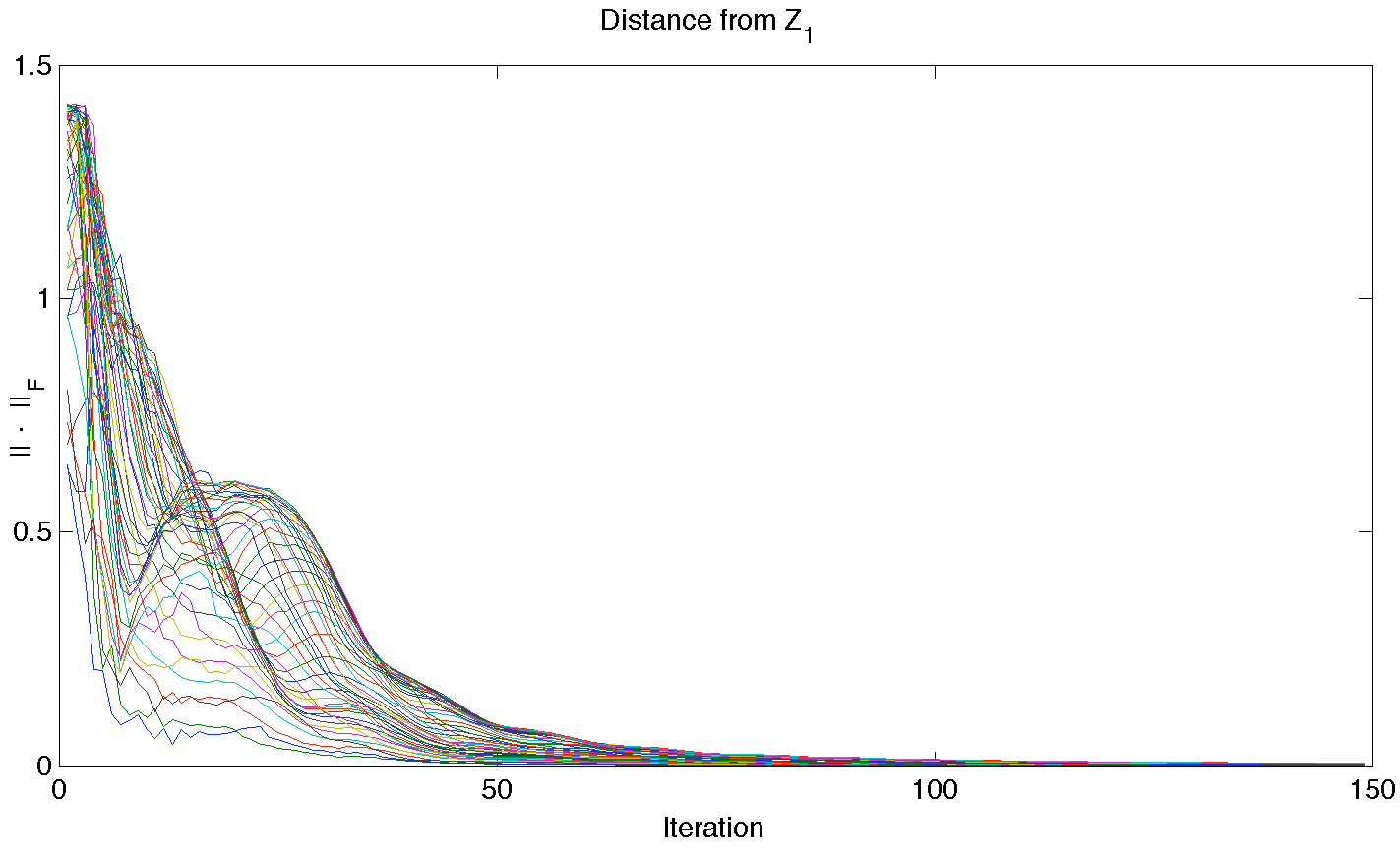}
\caption{\label{fig:admm_dec}
Error between average rotation estimates for a satellite formation of $N=50$ satellites with nearest- and second- neighbor communication using the distributed ADMM algorithm. The gradient weight is set as $\alpha=0.2$. 
}
\end{center}
\end{figure}

%

\begin{figure}
\begin{center}
\begin{tabular}{cc}
\includegraphics[scale=.4]{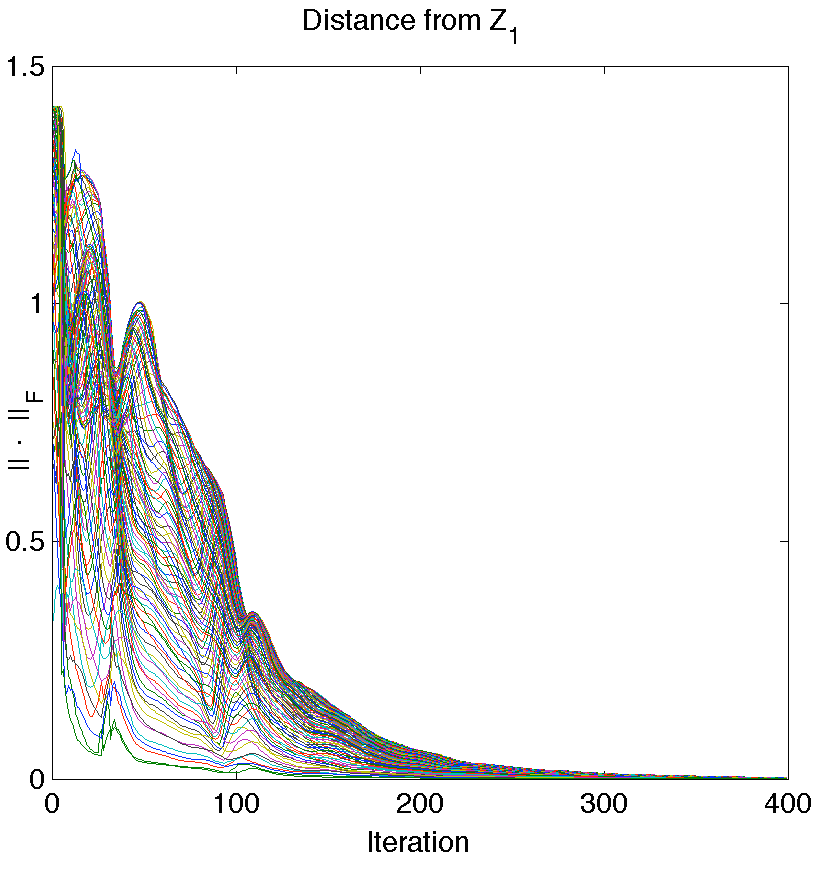} &
\includegraphics[scale=.4]{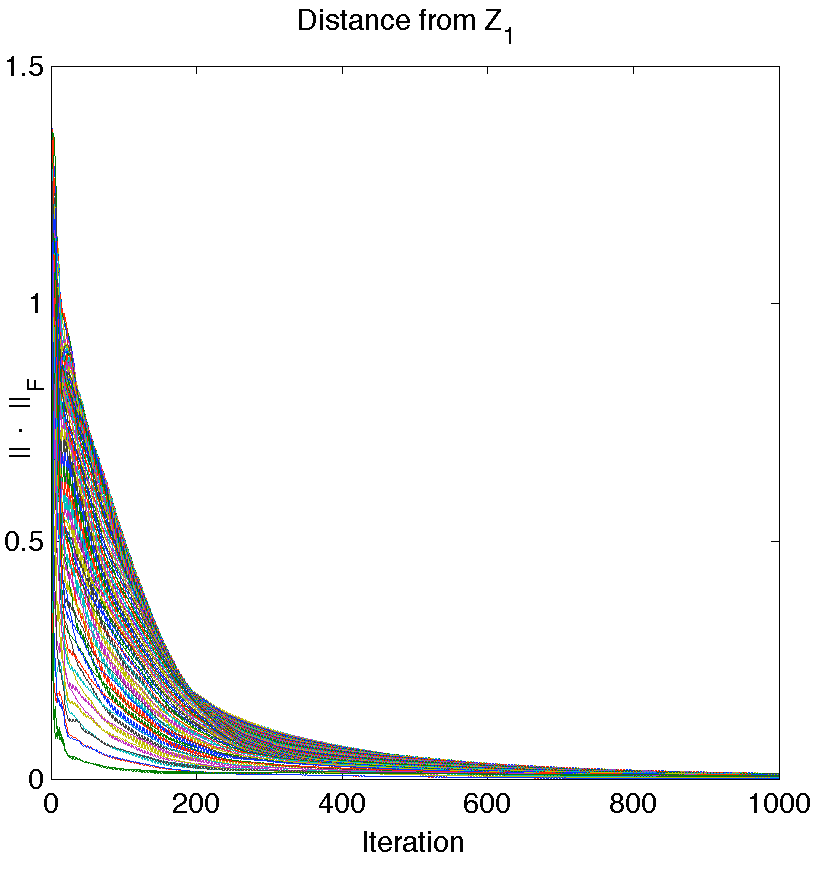} 
\end{tabular}
\caption{\label{fig:admm_dec_100}
Error between average rotation estimates for a satellite formation of $N=100$ satellites arranged in a ring topology with nearest neighbor (left) and second nearest neighbor (right) communication using the decentralized ADMM algorithm. The gradient descent weight is set as $\alpha=0.2$.
}
\end{center}
\end{figure}

The problem of synchronizing the attitude of $N=50$ satellites was then repeated using the distributed ADMM approach of \S \ref{sec:dist_admm}. Again, a ring topology was enforced, with Figure \ref{fig:admm_dec} showing the convergence of the algorithm when only allowing nearest neighbor and also second nearest neighbor communication. It is seen that the distributed ADMM solution is quicker and appears more regular than the dual ascent method in practice, although it is slightly more computationally demanding requiring both an eigen-decomposition and a projection onto the probability simplex.

Finally, to demonstrate the scalability properties of the algorithm, a $N=100$ satellite example with nearest- and second nearest-neighbor connectivity was simulated, with results shown in Figures \ref{fig:admm_dec_100}.

\subsection{Parallel Pose Estimation}
\begin{figure}
\begin{center}
\begin{tabular}{ccc}
\includegraphics[scale=0.6]{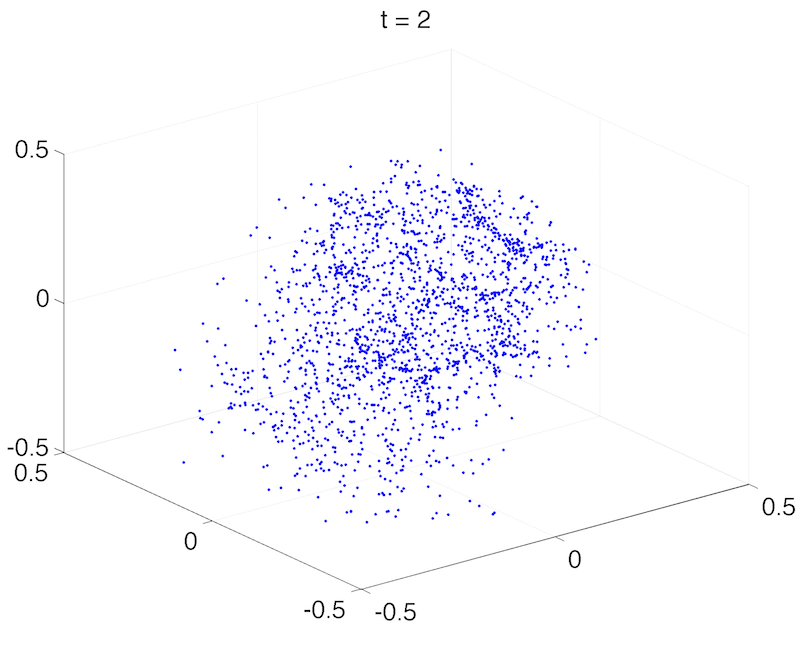} & 
\includegraphics[scale=0.6]{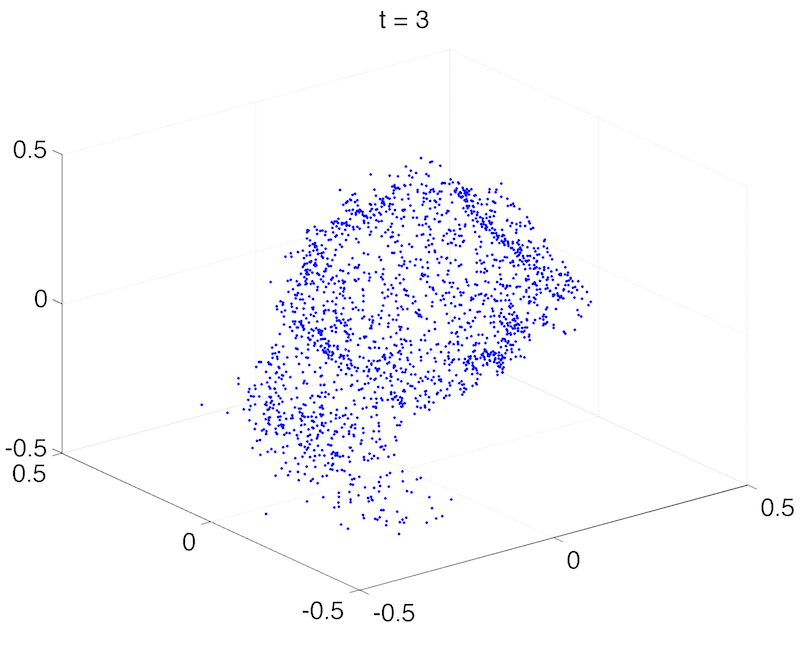} & 
\includegraphics[scale=0.6]{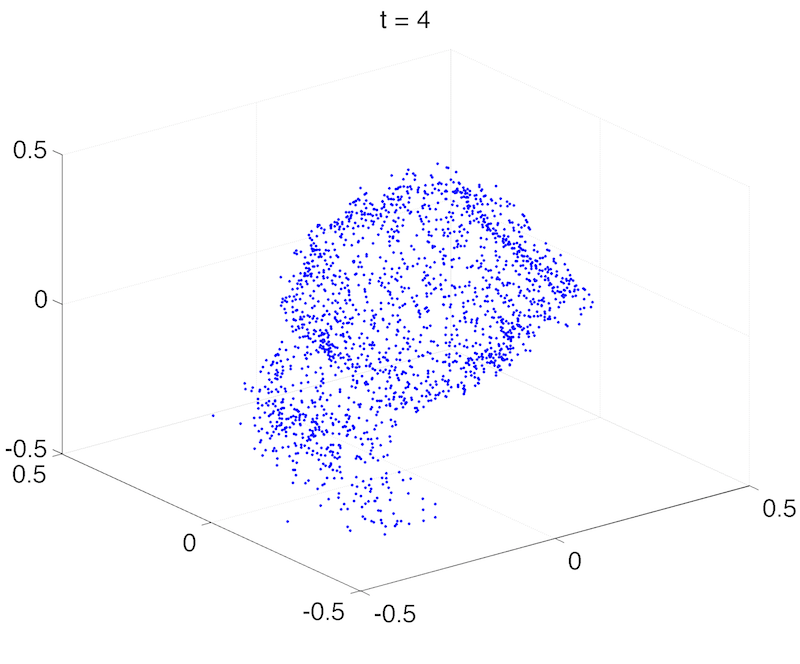} \\
\includegraphics[scale=0.6]{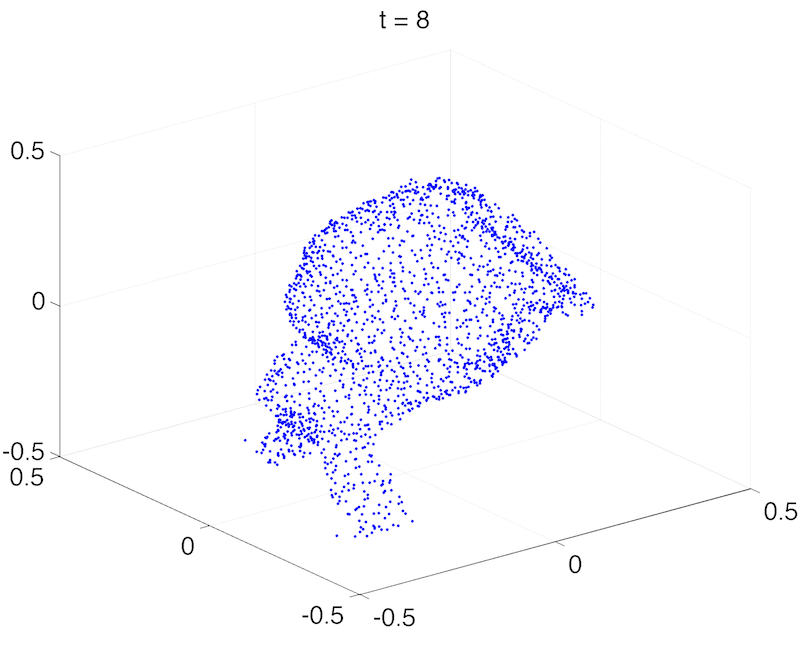} &
\includegraphics[scale=0.6]{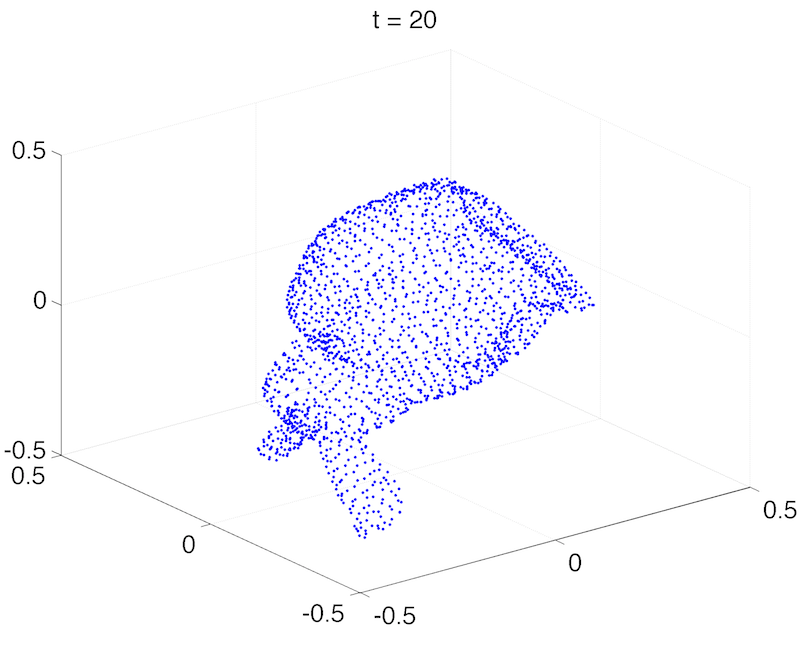} &
\includegraphics[scale=0.6]{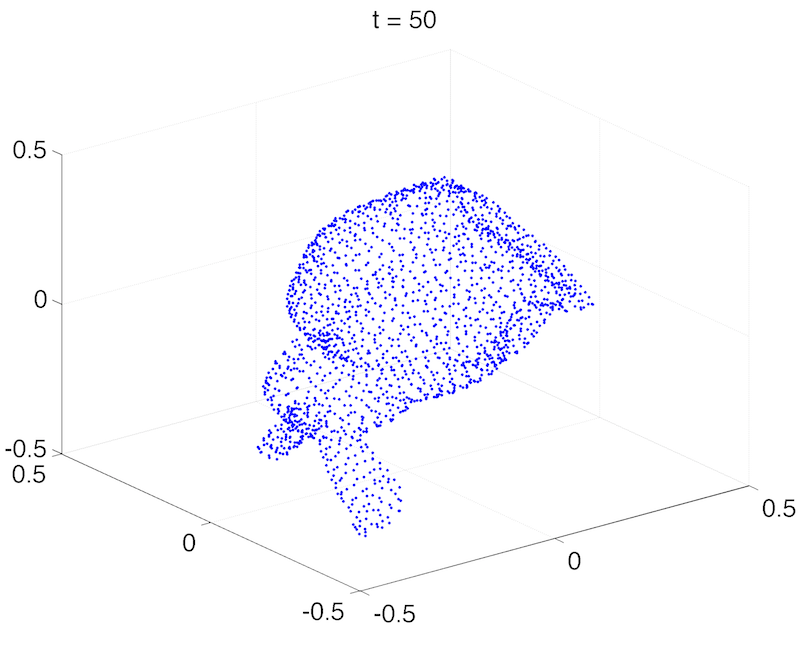}
\end{tabular}
\caption{\label{fig:bunny_lapse}
Visualization of the estimated Stanford Bunny model for $N=8$ different sensors. The inconsistency of the local estimates cause the local models (groups of vertices) to appear un-aligned. As consensus is achieved, the model becomes consistent.
}
\end{center}
\end{figure}

\begin{figure}
\begin{center}
\includegraphics[scale=1.0]{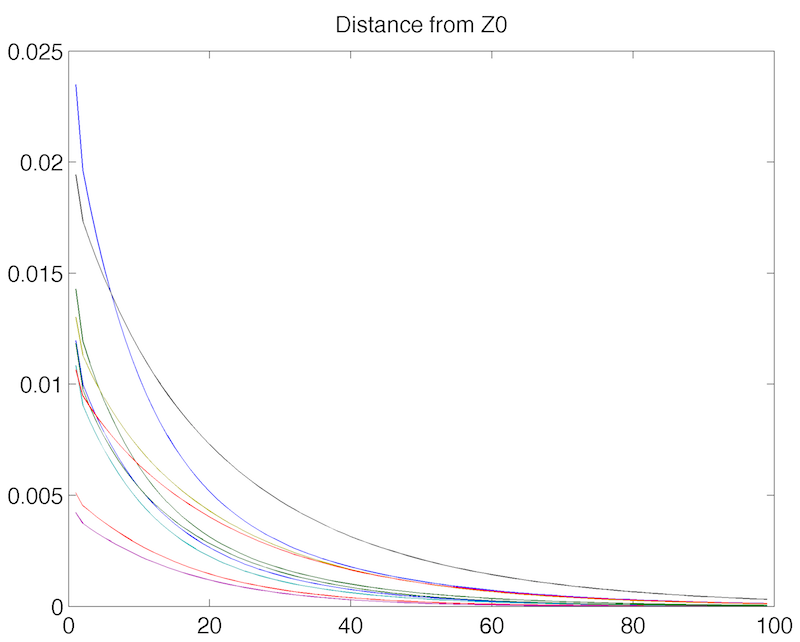}
\caption{\label{fig:admm_convergence}
Convergence of Algorithm \eqref{eq:admm} for $n=8$ disparate sensors for the Stanford Bunny model.
}
\end{center}
\end{figure}

The semi-distributed ADMM scheme is applied to the problem of pose estimation from computer vision. This problem is to identify a transformation that, when applied to a known object model, yields the object as perceived through the system's sensors. The convex relaxation has been applied previously in this context \cite{horowitz2014convex}. However, the approach presented in \S\ref{sec:admm} provides computationally appealing benefits, and arises in a number of ways. The first is to obtain a consistent estimate from a number of disparate sensors, for instance from a collection of cameras viewing a common object. This may be further extended for use in sensor fusion \cite{xiao2005scheme}, wherein the disparate sensors may provide different estimation modalities, and weights of the data matrices $B^i$, $C^i$ in \eqref{eq:non_convex} may correspond to priors on data strength. Alternatively, this consensus problem may arise when the pose estimation problem is split in order to take advantage of parallel processing capabilities, as follows. 

Pose estimation, i.e. Example 2 of \S\ref{sec:problem}, is solved for a 3D model of the Stanford Bunny, which has $s=1889$ vertices, and which is normalized to have unit dimensions. A random orientation is applied to the model, and the data is corrupted by Gaussian noise with standard deviation $\sigma=0.05$. The data is split between $N=10$ processing nodes and solved using Algorithm \eqref{eq:admm} where the dual variable weighting is set to $\alpha=2$. 

A visualization of the consensus progress is shown in Figure \ref{fig:bunny_lapse}, where it is seen that the differing estimates of the rotation cause the estimated model points to appear inconsistent. As a common estimate is obtained, it is seen that the points indeed do consistently match the model. Figure \ref{fig:admm_convergence} demonstrates the convergence of the error measure.

\begin{figure}
\begin{center}
\begin{tabular}{cc}
\includegraphics[scale=.4]{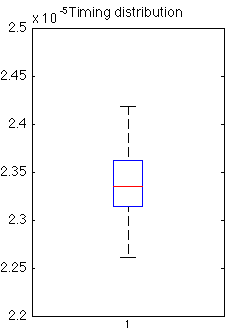} &
\includegraphics[scale=.4]{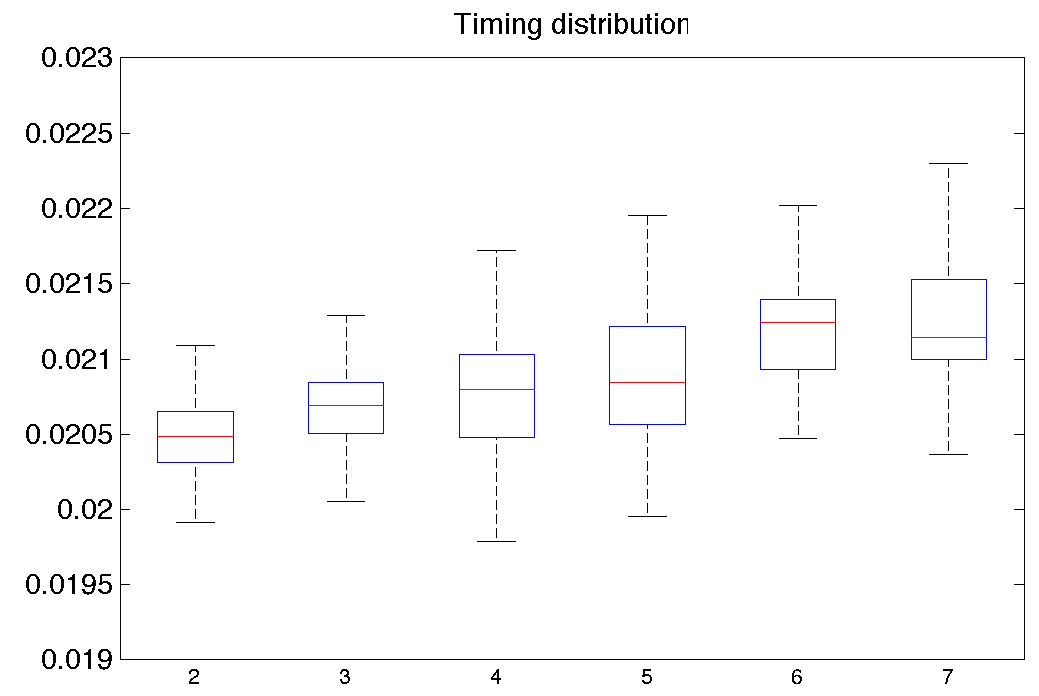} 
\end{tabular}
\caption{\label{fig:admm_subset_timing}
Timing results when varying the number of active cameras for the pose estimation problem.  On each box, the central mark is the median, the edges of the box are the 25th and 75th percentiles, the whiskers extend to the most extreme data points not considered outliers, and outliers are plotted individually.
}
\end{center}
\end{figure}

\begin{figure}
\begin{center}
\includegraphics[scale=.4]{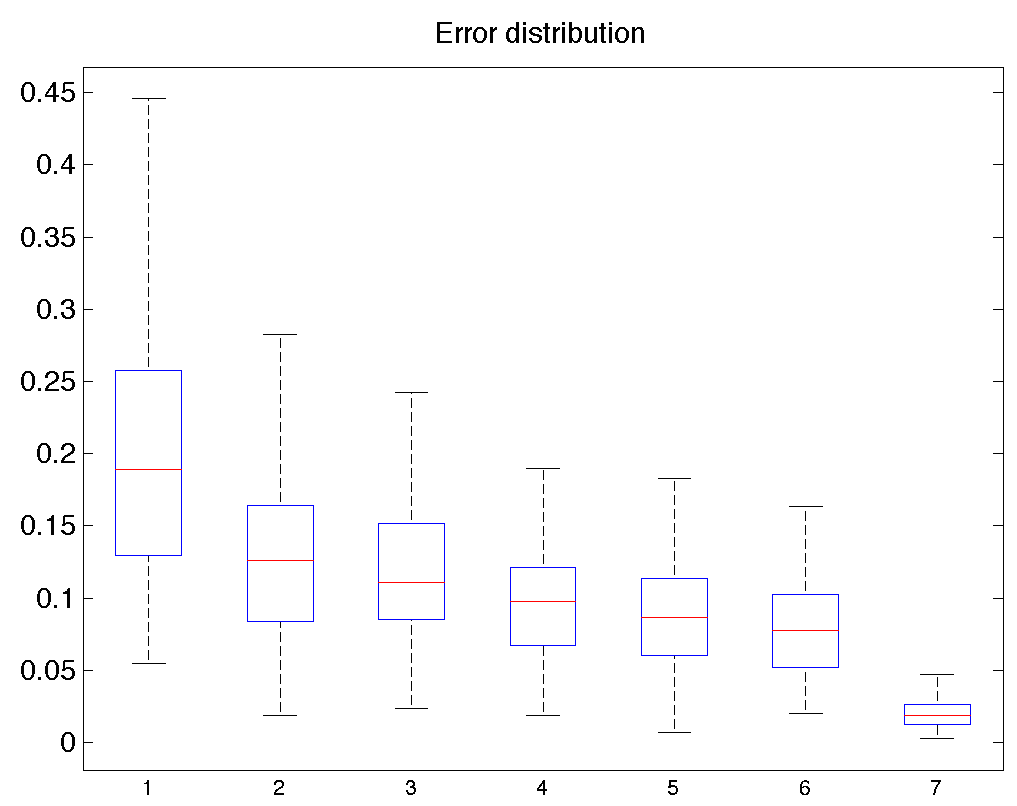}
\caption{\label{fig:admm_subset_error}
Estimation error results when varying the number of active cameras for the pose estimation problem.  On each box, the central mark is the median, the edges of the box are the 25th and 75th percentiles, the whiskers extend to the most extreme data points not considered outliers, and outliers are plotted individually.
}
\end{center}
\end{figure}

The problem was then examined by varying the number of virtual cameras available to collect observation data. The model was down sampled to $s=728$ vertices, and the standard deviation of the corrupting noise increased to $\sigma=0.10$. The observed data was then divided among $N=7$ sensors, and the estimation process repeated for $t=100$ trials for varying numbers of these sensors being active. The error and timing results are shown in Figures \ref{fig:admm_subset_timing} \& \ref{fig:admm_subset_error}.  We again emphasize that the timing results are for unoptimized, interpreted (Matlab) code running in \emph{serial} -- significant speed up is expected with compiled code and parallelized implementations.

\begin{rem}
Note that for the single camera case, the problem is no longer distributed. The computation of Lemma \ref{lem:analytic} is instead used, explaining the significant difference in computational time.
\end{rem}

\subsection{$SO(6)$ Consensus}
\begin{figure}
\begin{center}
\includegraphics[scale=.4]{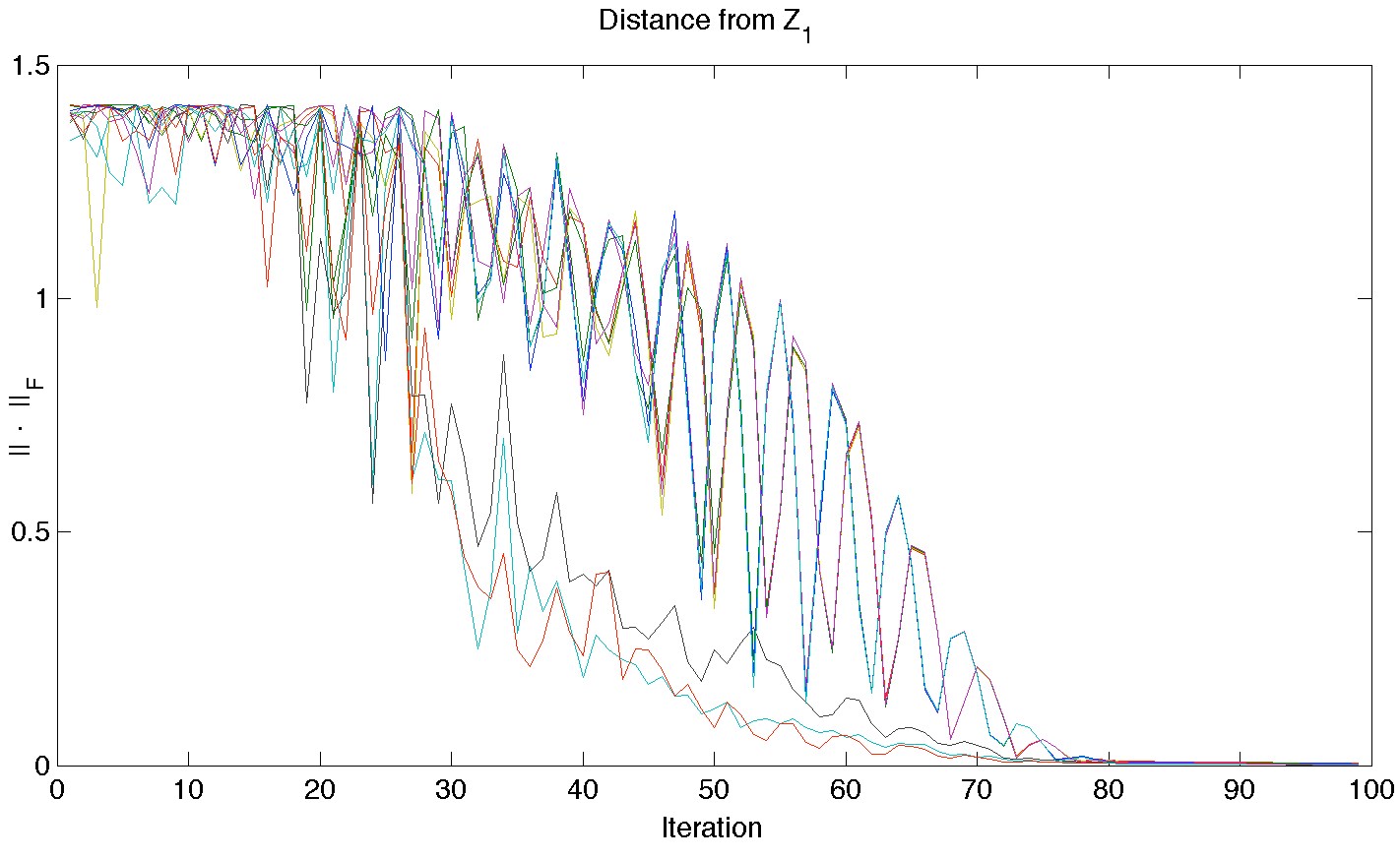}
\caption{\label{fig:so6_dual_decomp}
Convergence of Algorithm \eqref{eq:dualdecomp} in calculating the rotational average for $n=6$ and number of nodes $N=12$. The gradient weight was set as $\alpha=2$.
}
\end{center}
\end{figure}

Although applications of the Special Orthogonal group are most common for $SO(2)$ and $SO(3)$, efficient characterizations of $SO(n)$ for $n>3$ are necessary for computation on high dimensional datasets \cite{thakur2008}. In this application, rigid body transformations are used to change the perspective of a projection onto a lower dimensional structure. With this in mind, we apply Algorithm \eqref{eq:dualdecomp} to the rotation averaging problem where $n=6$, and number of nodes $N=12$, with the results shown in Figure \ref{fig:so6_dual_decomp}.

\section{Conclusion}
\label{sec:conclusion}
In this paper we have presented a convex approach to consensus on $SO(n)$ by leveraging recent results on the spectrahedral representation of $\co{SO(n)}$ and techniques from distributed optimization.  In particular, we show how meaningful consensus problems can be formulated as optimizations of linear functionals over the convex hull of the group of rotation matrices, and that these optimizations can be solved in a  distributed manner using dual decomposition and ADMM techniques -- these distributed computation methods are then used to implement a consensus protocol among the agents.  The benefits of our approach include rapid and guaranteed convergence to globally optimal consensus points, scalability and mathematical transparency.  Future work will focus on understanding how the choice of $\alpha$ in the ADMM algorithm affects boundary solution guarantees.


\bibliographystyle{IEEEtran}
\bibliography{biblio/comms,biblio/decentralized,biblio/matni,biblio/convex_vision,biblio/regularization,biblio/sonconsensus}

\end{document}